\documentclass[11pt]{article}
\usepackage[a4paper,margin=2.55cm]{geometry}
\usepackage{amsmath,amssymb,amsthm,mathtools}
\usepackage{microtype}
\usepackage{enumitem}
\usepackage{booktabs}
\usepackage{tabularx}
\usepackage{graphicx}
\usepackage{xcolor}
\usepackage[colorlinks=true,linkcolor=blue!55!black,citecolor=blue!55!black,urlcolor=blue!55!black,bookmarks=true,bookmarksnumbered=true]{hyperref}
\newtheorem{theorem}{Theorem}[section]
\newtheorem{lemma}[theorem]{Lemma}
\newtheorem{proposition}[theorem]{Proposition}
\newtheorem{corollary}[theorem]{Corollary}
\theoremstyle{definition}

\newtheorem{remark}[theorem]{Remark}
\DeclareMathOperator{\prox}{prox}
\DeclareMathOperator{\argmin}{argmin}
\newcommand{\Hil}{\mathcal H}
\newcommand{\R}{\mathbb R}
\newcommand{\eps}{\varepsilon}
\newcommand{\kap}{\kappa}
\title{Strong Hierarchical Selection for Inertial Forward--Backward at the Critical Tikhonov Exponent}
\author{Yekini Shehu\thanks{\raggedright School of Mathematical Sciences, Zhejiang Normal University, Jinhua 321004, China; e-mail: yekini.shehu@zjnu.edu.cn, deltanougt2006@yahoo.com}}
\date{}
\hypersetup{pdftitle={Strong Selection in the Critical Inertial Forward-Backward Algorithm},pdfauthor={Yekini Shehu},pdfsubject={Strong hierarchical selection at the critical Tikhonov exponent},pdfkeywords={inertial forward-backward, Tikhonov regularization, strong convergence, hierarchical selection}}
\begin{document}
\maketitle
\begin{abstract}
\noindent
We establish strong hierarchical selection and endpoint convergence rates for a Nesterov-type inertial forward--backward iteration at the critical Tikhonov exponent, with a single anchor $c k^{-2}(y_k-a)$ evaluated at the extrapolated point. For a general composite convex objective in a real Hilbert space, and for every $\alpha>3$, $c>0$, anchor $a$, and stepsize $0<s<1/L$, the iterates converge strongly to the projection of $a$ onto the solution set, with $F(x_k)-F^*=o(k^{-2})$ and $\|x_{k+1}-x_k\|=o(k^{-1})$. The conclusions hold for every admissible initialization, with no geometric restriction on the solution set or a priori boundedness assumption on the trajectory. The proof develops a critical-scale energy analysis: weighted square-summability of the Tikhonov-path increments controls the moving-target terms, while a sign-sensitive bootstrap reduces the generally nonsummable radial forcing to $o(k^{-1})$. A jointly coercive regularized energy and a two-index logarithmic convergence lemma then yield strong selection and both endpoint rates. The same analysis covers additive errors $e_k$ in the proximal arguments whenever $\sum_k k\|e_k\|<\infty$, preserving strong selection and both rates.
\end{abstract}



\section{Introduction}\label{sec:intro}
\subsection{Problem and algorithm}
Let $\Hil$ be a real Hilbert space with inner product $\langle\cdot,\cdot\rangle$ and induced norm $\|\cdot\|$. We study
\begin{equation}\label{eq:problem}
 \min_{x\in\Hil}F(x),\qquad F=f+g,
\end{equation}
where $f:\Hil\to\R$ is convex and differentiable with an $L$-Lipschitz gradient, and $g:\Hil\to(-\infty,+\infty]$ is proper, convex, and lower semicontinuous. We take $L>0$; if the optimal Lipschitz constant is zero, any positive $L$ may be used. The solution set $S:=\argmin F$ is nonempty and $F^*:=\min F$. For an anchor $a\in\Hil$, write $P_Sa$ for the unique projection of $a$ onto the closed convex set $S$.

The iteration is
\begin{equation}\label{eq:algorithm}
\boxed{
\begin{aligned}
 y_k&=x_k+\frac{k-1}{k+b}(x_k-x_{k-1}),\qquad b=\alpha-1,\\
 x_{k+1}&=\prox_{sg}\!\left(y_k-s\nabla f(y_k)-\frac{sc}{k^2}(y_k-a)\right),
\end{aligned}}
\end{equation}
with $x_0=x_1\in\operatorname{dom}g$ and $0<s<1/L$. Here $\prox_{sg}(v):=\argmin_{u\in\Hil}\{g(u)+(2s)^{-1}\|u-v\|^2\}$. Each iteration uses one gradient evaluation and one proximal evaluation.

The limit $P_Sa$ minimizes $\frac12\|x-a\|^2$ over $S$. For $a=0$, this is minimum-norm selection. For general $a$, it is selection relative to a prescribed state. This hierarchical principle has a long history in proximal and Tikhonov approximation; see, for example, \cite{Cabot2005,CominettiPeypouquetSorin2008,BotCsetnekMeier2019}.

\subsection{Position relative to prior work}
Inertial proximal methods arise from damped oscillator discretizations \cite{AlvarezAttouch2001}. Nesterov's method introduced the accelerated $O(k^{-2})$ objective rate \cite{Nesterov1983}, FISTA extended it to composite objectives \cite{BeckTeboulle2009}, and Su--Boyd--Cand\`es identified its vanishing-damping interpretation \cite{SuBoydCandes2016}. Iterate convergence was proved by Chambolle and Dossal \cite{ChambolleDossal2015}; Attouch and Peypouquet subsequently obtained the refined $o(k^{-2})$ objective rate and the corresponding discrete-velocity rate \cite{AttouchPeypouquet2016}.

The iteration belongs to the family of Nesterov-type forward--backward methods with vanishing Tikhonov regularization. Its closest general discrete predecessor, due to Attouch, Cabot, Chbani, and Riahi \cite{AttouchCabotChbaniRiahi2018}, includes
\begin{equation}\label{eq:accr-scheme}
 y_k=x_k+\alpha_k(x_k-x_{k-1}),\qquad
 x_{k+1}=\prox_{sg}\!\left(y_k-s\nabla f(y_k)-s\eps_ky_k\right).
\end{equation}
For $\alpha_k=(k-1)/(k+b)$ and $\eps_k=ck^{-2}$, this is precisely the zero-anchor case of \eqref{eq:algorithm}. The translation $\widetilde x_k=x_k-a$, $\widetilde f(w)=f(w+a)$, and $\widetilde g(w)=g(w+a)$ yields the arbitrary-anchor formulation from the zero-anchor iteration for $\widetilde F=\widetilde f+\widetilde g$, whose minimum-norm solution is $P_Sa-a$. This formulation expresses hierarchical selection relative to any prescribed state $a$ and recovers minimum-norm selection at $a=0$.

The issue addressed here is the \emph{last-iterate theory at the critical exponent} $p=2$ for $\eps_k=ck^{-p}$. At this exponent, $k^2\eps_k$ remains of order one: anchoring and accelerated objective decay act on the same energy scale. The Tikhonov term therefore cannot be treated as a perturbation that simply disappears from the $k^{-2}$ analysis.

This boundary is visible in the closest discrete results. The Tikhonov application in \cite{AttouchCabotChbaniRiahi2018} obtains strong ergodic minimum-norm convergence for coefficients that do not vanish too rapidly, rather than the full critical last-iterate conclusion proved here. Attouch and L\'aszl\'o \cite{AttouchLaszlo2024} obtained rates $o(k^{-2s})$ and $o(k^{-s})$ for every $0<s<1$ in critical proximal regimes, together with a general $\liminf$ selection statement; full strong selection was established under an eventual one-sided ball alternative. L\'aszl\'o's inertial proximal analysis \cite{LaszloInertialProximal2025} likewise leaves full critical selection open in the relevant single-anchor regime. His proximal--gradient method \cite{LaszloProxGradient2025} uses two Tikhonov terms and obtains $O(k^{-2})$ at critical parameters, but no critical iterate-convergence conclusion. Adly, Ho, and Nguyen \cite{AdlyHoNguyen2026} proved strong minimum-norm selection for an explicit smooth method when $0<p<2$; at $p=2$, their theorem retains $O(k^{-2})$ objective and $O(k^{-1})$ velocity rates without guaranteeing strong selection. Double-anchor methods, such as that of Karapetyants and L\'aszl\'o \cite{KarapetyantsLaszlo2024}, belong to a different noncritical regime.

Related continuous-time analyses include fixed-damping and vanishing-damping Tikhonov dynamics \cite{AttouchChbaniRiahi2018,AttouchBalhagChbaniRiahi2022}, Hessian-driven and Newton-like models \cite{BotCsetnekLaszlo2021,BotCsetnekLaszlo2024}, and monotone-inclusion frameworks \cite{BotGradMeierStaudigl2021,BotNguyen2026Flows}. At the critical continuous scale, Li and Xiao recently proved strong minimum-norm selection without eventual-ball conditions, with $o(t^{-2})$ residual and $o(t^{-1})$ velocity rates, for mixed-order primal--dual dynamics \cite{LiXiao2026Mixed}; they also obtained critical strong-selection and fast-rate results for nonlinear saddle-point dynamics and an associated inexact algorithm \cite{LiXiao2026Saddle}. Those results concern different dynamics and discretizations. They show that the critical little-$o$ endpoint is a current and substantive analytical question, while leaving the explicit composite iteration \eqref{eq:algorithm} as the object of the present proof. For completeness, unanchored critical Nesterov dynamics have separate iterate-convergence results in \cite{BotFadiliNguyen2025,JangRyu2026}; those results select no anchor-dependent limit. Different anchored bilevel and monotone-equation schemes are analyzed in \cite{BotChencheneCsetnekHulett2026,BotChenchene2026Anchoring}.

For the single-anchor composite iteration \eqref{eq:algorithm} at $p=2$, Theorem~\ref{thm:main} establishes full last-iterate strong selection together with the endpoint rates $o(k^{-2})$ and $o(k^{-1})$. The critical-scale analysis introduces recurrence estimates that control the moving-path terms and radial forcing while avoiding both a geometric alternative and an a priori boundedness assumption.

\subsection{Main result}
\begin{theorem}[Critical strong hierarchical selection]\label{thm:main}
Let $\alpha>3$, $c>0$, $a\in\Hil$, and $0<s<1/L$. Then the sequence generated by \eqref{eq:algorithm} satisfies
\begin{equation}\label{eq:main-conclusion}
 x_k\longrightarrow P_Sa\quad\text{strongly},\qquad
 F(x_k)-F^*=o(k^{-2}),\qquad
 \|x_{k+1}-x_k\|=o(k^{-1}).
\end{equation}
The same conclusions hold for every initialization $x_0=x_1\in\operatorname{dom}g$, with no assumption beyond those stated for \eqref{eq:problem}.
\end{theorem}

Theorem~\ref{thm:inexact} extends the same conclusions to additive proximal-argument errors satisfying $\sum_k k\|e_k\|<\infty$. This error analysis uses the same critical-regime mechanism rather than a separate boundedness hypothesis.

\begin{remark}[The selection content in one dimension]
Take $\Hil=\R$, $g=0$, and
\[
 f(x):=\frac12\operatorname{dist}(x,[-1,1])^2.
\]
Then $L=1$, $S=[-1,1]$, and $P_Sa=\min\{1,\max\{-1,a\}\}$. If $a>1$, the regularized path introduced below is explicitly
\[
 \xi_k=\frac{1+\eps_k a}{1+\eps_k}
 =1+\frac{\eps_k}{1+\eps_k}(a-1)\longrightarrow1=P_Sa.
\]
Thus selection concerns a prescribed endpoint of the solution interval, not convergence to an unspecified minimizer. The general proof uses neither this explicit path nor finite-dimensionality.
\end{remark}

\subsection{The analytical obstruction and the proof mechanism}\label{subsec:mechanism}
The proof can be understood without the detailed coefficients. Define the regularized objectives and their minimizers by
\[
 \Phi_k=F+\frac{\eps_k}{2}\|\cdot-a\|^2,\qquad \xi_k=\argmin\Phi_k,
 \qquad \eps_k=ck^{-2}.
\]
A standard inertial energy can be adapted to a \emph{fixed} target and a \emph{fixed} objective. At iteration $k$, however, both the target $\xi_k$ and the objective $\Phi_k$ move. After the exact index shifts and absorption of the path increments, the analysis yields a recurrence of the form
\begin{equation}\label{eq:schematic-recurrence}
 H_{k+1}\le H_k-\frac{\rho}{k}D_k
 +\frac{C_q}{k}q_{k+1,+}+\zeta_k,
 \qquad H_k\le C_0(D_k+D_{k-1}),
\end{equation}
where $H_k$ is a coercive regularized energy, $D_k$ combines objective-gap and velocity dissipation, $\sum_k\zeta_k<\infty$, and
\[
 q_k:=\|\xi_k-a\|^2-\|x_k-a\|^2,\qquad q_{k,+}:=\max\{q_k,0\}.
\]
The precise statement is Proposition~\ref{prop:absorption}. Three points explain why the critical case does not follow from a routine summability argument.

First, replacing $\xi_k$ by $\xi_{k+1}$ creates inner products with $h_k:=\xi_{k+1}-\xi_k$. Strong convexity of successive regularized objectives gives the automatic estimate
\begin{equation}\label{eq:intro-path-sum}
 \sum_{k\ge1}k\|h_k\|^2
 \le\|P_Sa-a\|^2-\|\xi_1-a\|^2<\infty.
\end{equation}
This weighted \emph{square}-summability is exactly what the absorbed energy recurrence requires; no unproved finite-path-length assumption $\sum_k\|h_k\|<\infty$ is made.

Second, changing $\Phi_k$ to $\Phi_{k+1}$ produces the radial forcing $C_qk^{-1}q_{k+1,+}$. The uniform bound $q_{k,+}\le\|P_Sa-a\|^2$ gives only an $O(k^{-1})$ forcing, which is generally nonsummable. It is nevertheless sufficient for a preliminary bootstrap: scaling $H_k$ by $k^{-2}$ gives $H_k/k^2\to0$, hence regularized-gap convergence and then $F(x_k)\to F^*$. Objective convergence now has a geometric consequence. A minimizing sequence cannot have a subsequence whose squared distance to $a$ remains strictly below $\|P_Sa-a\|^2$: boundedness of that subsequence and weak lower semicontinuity would produce a minimizer closer to $a$ than its projection. Consequently $q_{k,+}\to0$, and the forcing improves from $O(k^{-1})$ to $o(k^{-1})$. This upgrade uses objective convergence, not the strong convergence that remains to be proved.

Third, $o(k^{-1})$ forcing is still not summable. The two-index domination in \eqref{eq:schematic-recurrence} supplies the missing principle. On any interval where $H_k$ stays above a fixed level, the pair $D_{k-1},D_k$ generates harmonic dissipation. A short logarithmic excursion cannot acquire a fixed positive gain from vanishing forcing, while a long logarithmic excursion incurs an arbitrarily large energy loss. Lemma~\ref{lem:two-index} formalizes this dichotomy and proves $H_k\to0$. Strong convergence and both little-$o$ rates then follow from joint coercivity of $H_k$ and strong convergence of the path $\xi_k$.

There is one further technical point. The chosen energy contains signed inertial coefficients, so positivity is not visible term by term. The objective gap and velocity terms instead form a jointly coercive quadratic form. Its limiting determinant is
\[
 \frac{b-\beta}{\beta}+\frac{sc}{2\beta^2}>0,
 \qquad 2<\beta<b,
\]
which explains both the free-inertia split and why no restriction on the anchor strength $c$ is needed.

These four ingredients--weighted path summability, the radial bootstrap, two-index logarithmic convergence, and joint coercivity--close the critical last-iterate endpoint argument. Together they yield strong hierarchical selection, both little-$o$ rates, and the inexact extension within a unified energy analysis.

\subsection{Organization and notation}
Section~\ref{sec:path} develops the regularized path and the radial lemma. Section~\ref{sec:discrete-lemmas} proves the two scalar convergence principles. Section~\ref{sec:fixed} derives the fixed-objective, fixed-target recurrence and establishes joint coercivity. Section~\ref{sec:moving} computes the moving-path shifts exactly and obtains \eqref{eq:schematic-recurrence}. Section~\ref{sec:proof} assembles the bootstrap and proves Theorem~\ref{thm:main}. Section~\ref{sec:errors} treats additive errors, and Section~\ref{sec:numerics} gives deterministic illustrations.

Throughout, $(t)_+:=\max\{t,0\}$ and $(t)_-:=\max\{-t,0\}$ for $t\in\R$. We write $A_k\asymp k^r$ when positive lower and upper bounds of order $k^r$ hold for all sufficiently large $k$. Constants denoted by $C$ may change from line to line but depend only on the fixed parameters and the finite initial segment, never on an assumed trajectory bound. Standard Hilbert-space convex analysis is used in the form presented in \cite{BauschkeCombettes2017}.
\section{The discrete Tikhonov path and radial selection}\label{sec:path}
\subsection{Regularized objectives and their minimizers}
Throughout, set $z:=P_Sa$ and $R_*:=\|z-a\|$. For $k\ge1$, define
\begin{equation}\label{eq:path-definition}
 \eps_k:=\frac{c}{k^2},\qquad
 \Phi_k(u):=F(u)+\frac{\eps_k}{2}\|u-a\|^2,\qquad
 \xi_k:=\argmin_{u\in\Hil}\Phi_k(u).
\end{equation}
Each $\Phi_k$ is proper, lower semicontinuous, and $\eps_k$-strongly convex. Coercivity and strong convexity therefore guarantee the existence and uniqueness of $\xi_k$. We write
\begin{equation}\label{eq:gap-radius}
 P_k:=\Phi_k(x_k)-\Phi_k(\xi_k),\qquad R_k:=\|\xi_k-a\|.
\end{equation}
The regularized gap is nonnegative, and strong convexity gives the quantitative lower bound
\begin{equation}\label{eq:strong-gap}
 P_k\ge\frac{\eps_k}{2}\|x_k-\xi_k\|^2.
\end{equation}
This estimate will be used in both directions: it turns gap dissipation into distance dissipation, and it converts cross terms involving $x_k-\xi_k$ into terms of order $kP_k$.

Comparing $\Phi_k(\xi_k)$ with $\Phi_k(z)$ gives
\[
 F(\xi_k)+\frac{\eps_k}{2}R_k^2
 \le F^*+\frac{\eps_k}{2}R_*^2.
\]
Since $F(\xi_k)\ge F^*$, we obtain
\begin{equation}\label{eq:path-comparison}
 F(\xi_k)\le F^*+\frac{\eps_k}{2}R_*^2,
 \qquad R_k\le R_*.
\end{equation}
Thus the path is bounded and minimizing, with no assumption on the algorithmic trajectory.

\begin{lemma}[Strong convergence of the path]\label{lem:path-convergence}
The path satisfies $\xi_k\to z$ strongly and $R_k\to R_*$.
\end{lemma}
\begin{proof}
By \eqref{eq:path-comparison}, $(\xi_k)$ is bounded and $F(\xi_k)\to F^*$. Let $\bar\xi$ be any weak sequential cluster point. Weak lower semicontinuity of $F$ gives $\bar\xi\in S$. The same comparison, together with weak lower semicontinuity of the norm, yields $\|\bar\xi-a\|\le R_*$. Since $z$ is the unique nearest point to $a$ in $S$, $\bar\xi=z$. Every weak sequential cluster point therefore equals $z$, and boundedness implies $\xi_k\rightharpoonup z$.

Weak lower semicontinuity and the upper bound in \eqref{eq:path-comparison} now give
\[
 R_*\le\liminf_kR_k\le\limsup_kR_k\le R_*.
\]
Hence $R_k\to R_*=\|z-a\|$. In a Hilbert space, weak convergence together with convergence of $\|\xi_k-a\|$ to $\|z-a\|$ implies strong convergence of $\xi_k-a$ to $z-a$, and therefore $\xi_k\to z$ strongly.
\end{proof}

\subsection{An automatic path estimate}
The next lemma is specific to the discrete Tikhonov path, not to the algorithm. Its importance is that the weighted squared increments telescope against the monotone radii $R_k^2$.

\begin{lemma}[Automatic path summability]\label{lem:path-summability}
For every $k\ge1$, with $h_k:=\xi_{k+1}-\xi_k$,
\begin{equation}\label{eq:path-one-step}
 R_{k+1}^2-R_k^2
 \ge\frac{\eps_k+\eps_{k+1}}{\eps_k-\eps_{k+1}}\|h_k\|^2.
\end{equation}
Consequently, $(R_k^2)$ is nondecreasing and
\begin{equation}\label{eq:path-sum}
 k\|h_k\|^2\le R_{k+1}^2-R_k^2,
 \qquad
 \sum_{k=1}^{\infty}k\|h_k\|^2\le R_*^2-R_1^2<\infty.
\end{equation}
\end{lemma}
\begin{proof}
Strong convexity of the successive objectives at their respective minimizers gives
\[
 \Phi_k(\xi_{k+1})-\Phi_k(\xi_k)\ge\frac{\eps_k}{2}\|h_k\|^2,
 \qquad
 \Phi_{k+1}(\xi_k)-\Phi_{k+1}(\xi_{k+1})\ge\frac{\eps_{k+1}}{2}\|h_k\|^2.
\]
Add these inequalities. The $F$-values cancel, leaving
\[
 \frac{\eps_k}{2}\bigl(R_{k+1}^2-R_k^2\bigr)
 +\frac{\eps_{k+1}}{2}\bigl(R_k^2-R_{k+1}^2\bigr)
 \ge\frac{\eps_k+\eps_{k+1}}2\|h_k\|^2,
\]
which is exactly \eqref{eq:path-one-step}. Monotonicity follows because the fraction on the right is positive.

For $\eps_j=cj^{-2}$, direct simplification gives
\[
 \frac{\eps_k-\eps_{k+1}}{\eps_k+\eps_{k+1}}
 =\frac{2k+1}{2k^2+2k+1},
 \qquad
 k\frac{2k+1}{2k^2+2k+1}\le1.
\]
Taking reciprocals in \eqref{eq:path-one-step} therefore yields $k\|h_k\|^2\le R_{k+1}^2-R_k^2$. Summing from $1$ to $K$ and using $R_{K+1}\le R_*$ proves \eqref{eq:path-sum}, and letting $K\to\infty$ completes the argument.
\end{proof}

\begin{remark}[What the path estimate does and does not assume]
The bound \eqref{eq:path-sum} is a weighted square-summability statement. It does not assert that $\sum_k\|h_k\|<\infty$, and no such finite-path-length property is needed later. It also does not use the iterates $(x_k)$: its constant depends only on the anchor-selected radius $R_*$ and the initial path point $\xi_1$.
\end{remark}

\subsection{The negative radial part}
The moving-objective calculation will produce the signed radial quantity
\begin{equation}\label{eq:q-definition}
 q_k:=R_k^2-\|x_k-a\|^2.
\end{equation}
Only $q_{k,+}$ has the unfavorable sign in the energy recurrence. The next lemma shows that objective convergence alone makes this positive part asymptotically negligible. It does not require boundedness of the full sequence $(x_k)$.

\begin{lemma}[Negative radial part]\label{lem:negative-radial}
If $F(x_k)\to F^*$, then
\begin{equation}\label{eq:negative-radial}
 \bigl(\|x_k-a\|^2-R_*^2\bigr)_-\longrightarrow0.
\end{equation}
\end{lemma}
\begin{proof}
Suppose the assertion fails. Then there are $\eta>0$ and a subsequence $(x_{k_j})$ such that
\[
 \|x_{k_j}-a\|^2\le R_*^2-\eta.
\]
This subsequence lies in a fixed closed bounded ball. After extracting a further subsequence, it converges weakly to some $\bar x$ in that ball. Since $F(x_{k_j})\to F^*$ and $F$ is weakly lower semicontinuous, $\bar x\in S$. Weak lower semicontinuity of the norm also gives $\|\bar x-a\|^2\le R_*^2-\eta$, contradicting the definition of $R_*=\min_{u\in S}\|u-a\|$.
\end{proof}

The connection with $q_k$ is immediate. Since $R_k\le R_*$, positivity of $q_k$ implies $\|x_k-a\|<R_*$ and
\[
 0<q_k=R_k^2-\|x_k-a\|^2
 \le R_*^2-\|x_k-a\|^2.
\]
Consequently,
\begin{equation}\label{eq:q-positive-bound}
 q_{k,+}\le\bigl(\|x_k-a\|^2-R_*^2\bigr)_-
 \quad\text{and always}\quad q_{k,+}\le R_*^2.
\end{equation}
The two bounds play different roles. The uniform second bound is enough for the initial scaled-energy bootstrap. After objective convergence has been proved, Lemma~\ref{lem:negative-radial} and the first bound give $q_{k,+}\to0$, upgrading the radial term $q_{k,+}/k$ to $o(k^{-1})$. This separation avoids a circular appeal to strong convergence.
\section{Two discrete convergence lemmas}\label{sec:discrete-lemmas}
This section isolates the two scalar principles used in the bootstrap. Both are independent of optimization. The first handles a damped recursion with summable forcing; the second handles a nonsummable $o(k^{-1})$ forcing and a two-index dissipation block.

\subsection{The scaled-energy bootstrap}
\begin{lemma}[Damped recursion with summable forcing]\label{lem:damped}
Let $E_k\ge0$ satisfy, for all sufficiently large $k$,
\begin{equation}\label{eq:damped-recursion}
 E_{k+1}\le\left(1-\frac{\alpha}{k}+\sigma_k\right)E_k+\eta_k,
\end{equation}
where $\alpha>0$, $\sigma_k\ge0$, $\sum_k\sigma_k<\infty$, $\eta_k\ge0$, and $\sum_k\eta_k<\infty$. Then $E_k\to0$.
\end{lemma}
\begin{proof}
Choose $m$ so large that the factors in \eqref{eq:damped-recursion} are nonnegative and
$1-\alpha/k+\sigma_k\le\exp(-\alpha/k+\sigma_k)$ for $k\ge m$. Unrolling the recursion and using $\sum_k\sigma_k<\infty$ give a constant $C_m$, independent of $k$, such that
\[
 E_k\le C_m\left(\frac{m}{k}\right)^\alpha E_m
 +C_m\sum_{j=m}^{k-1}\left(\frac{j}{k}\right)^\alpha\eta_j.
\]
For fixed $N>m$, split the last sum at $N$. The bound $C_mk^{-\alpha}\sum_{j=m}^{N-1}j^\alpha\eta_j$ on the early part tends to zero. The late part is at most $C_m\sum_{j=N}^{\infty}\eta_j$. First let $k\to\infty$ and then $N\to\infty$. Summability of $(\eta_j)$ makes the second bound arbitrarily small.
\end{proof}

In the application, $E_k=H_k/k^2$ and $\alpha=2$. The factor $k^2/(k+1)^2$ differs from $1-2/k$ by a nonnegative $O(k^{-2})$ correction, which is absorbed into $\sigma_k$. Thus no exact product formula is required.

\subsection{A logarithmic principle for two-index dissipation}
The second lemma is deliberately formulated for an energy dominated by two consecutive dissipation terms. The input $r_k=o(k^{-1})$ is not assumed summable, so direct summation cannot prove convergence.

\begin{lemma}[Two-index logarithmic convergence]\label{lem:two-index}
Let $H_k,D_k\ge0$ and suppose that, for all sufficiently large $k$,
\begin{equation}\label{eq:two-index-hypotheses}
 H_{k+1}\le H_k-\frac{\rho}{k}D_k+r_k+\zeta_k,
 \qquad H_k\le C(D_k+D_{k-1}),
\end{equation}
where $\rho,C>0$, $r_k\ge0$, $r_k=o(k^{-1})$, $\zeta_k\ge0$, and $\sum_k\zeta_k<\infty$. Then $H_k\to0$.
\end{lemma}
\begin{proof}
Write $w_k:=r_k+\zeta_k$, so that $w_k\to0$. Assume, for contradiction, that $\limsup_kH_k>4\eps$ for some $\eps>0$. There are two possibilities.

\medskip
\noindent\emph{Case 1: $H_k\ge2\eps$ eventually.}
For every large $k$, the two-index domination gives
\begin{equation}\label{eq:pair-dissipation}
 \frac{D_{k-1}}{k-1}+\frac{D_k}{k}
 \ge\frac{D_{k-1}+D_k}{k}
 \ge\frac{H_k}{Ck}\ge\frac{2\eps}{Ck}.
\end{equation}
Adding the recurrence at indices $k-1$ and $k$ therefore gives
\[
 H_{k+1}-H_{k-1}
 \le-\frac{2\rho\eps}{Ck}+o(k^{-1})+\zeta_{k-1}+\zeta_k.
\]
Sum from $k=m$ to $K$. The harmonic part contributes $-\frac{2\rho\eps}{C}\bigl(\ln(K/m)+O(1)\bigr)$, the forcing contributes $\sum_{k=m}^{K}o(k^{-1})=o(\ln K)$, and the summable part contributes $O(1)$; hence the right-hand side tends to $-\infty$ as $K\to\infty$. The telescoped left-hand side $H_{K+1}+H_K-H_m-H_{m-1}$ remains bounded below by $-(H_m+H_{m-1})$, since all $H_j$ are nonnegative. This is impossible.

\medskip
\noindent\emph{Case 2: $H_k<2\eps$ infinitely often.}
Choose times $n_j\to\infty$ with $H_{n_j}\ge4\eps$. For each large $j$, let $l_j<n_j$ be the last index before $n_j$ with $H_{l_j}<2\eps$, and put $m_j:=l_j+1$. Thus
\begin{equation}\label{eq:excursion-level}
 H_{m_j-1}<2\eps,
 \qquad H_k\ge2\eps\quad(m_j\le k\le n_j).
\end{equation}
Dropping the dissipation in the recursion gives the fixed positive excursion gain
\begin{equation}\label{eq:excursion-gain}
 2\eps<H_{n_j}-H_{m_j-1}
 \le\sum_{k=m_j-1}^{n_j-1}w_k.
\end{equation}
Define the logarithmic lengths
\[
 L_j:=\sum_{k=m_j-1}^{n_j-1}\frac1k.
\]
If $(L_j)$ had a bounded subsequence, then along that subsequence $m_j\to\infty$ and
\[
 \sum_{k=m_j-1}^{n_j-1}r_k
 \le L_j\sup_{k\ge m_j-1}(kr_k)\to0,
 \qquad
 \sum_{k=m_j-1}^{n_j-1}\zeta_k\to0.
\]
This contradicts \eqref{eq:excursion-gain}. Hence $L_j\to\infty$.

On the other hand, \eqref{eq:excursion-level} and the pairwise estimate \eqref{eq:pair-dissipation} yield
\[
 2\sum_{k=m_j-1}^{n_j-1}\frac{D_k}{k}
 \ge\sum_{k=m_j}^{n_j-1}
 \left(\frac{D_{k-1}}{k-1}+\frac{D_k}{k}\right)
 \ge\frac{2\eps}{C}\sum_{k=m_j}^{n_j-1}\frac1k.
\]
The first inequality counts each interior dissipation term at most twice and includes the two endpoint terms on the left. Summing the original recurrence from $m_j-1$ to $n_j-1$ consequently gives
\begin{equation}\label{eq:excursion-loss}
 H_{n_j}-H_{m_j-1}
 \le-\frac{\rho\eps}{C}\sum_{k=m_j}^{n_j-1}\frac1k
 +\sum_{k=m_j-1}^{n_j-1}r_k
 +\sum_{k=m_j-1}^{n_j-1}\zeta_k.
\end{equation}
The first sum on the right equals $L_j+O(1)$ and tends to $+\infty$ with $j$. For the second, fix $K_0$ such that $kr_k\le\rho\eps/(2C)$ for all $k>K_0$ and split the sum at $K_0$: the initial segment contributes $O(1)$ independently of $j$, while the tail is at most $\frac{\rho\eps}{2C}L_j$; this estimate covers both $m_j\to\infty$ and bounded $(m_j)$. The third sum is bounded by $\sum_{k\ge1}\zeta_k=O(1)$ (and tends to zero whenever $m_j\to\infty$). Therefore the right-hand side of \eqref{eq:excursion-loss} is at most $-\frac{\rho\eps}{2C}L_j+O(1)$, which tends to $-\infty$ since $L_j\to\infty$. But its left-hand side is greater than $2\eps$ by \eqref{eq:excursion-gain}, a contradiction.

Neither case is compatible with $\limsup_kH_k>4\eps$. Since $\eps>0$ was arbitrary and $H_k\ge0$, $H_k\to0$.
\end{proof}

\begin{remark}[Why the logarithmic argument is necessary]
The recursion alone permits long intervals on which $H_k$ remains above a fixed level. On every such interval, however, the two-index domination converts that level into dissipation of order $1/k$. A high-level excursion therefore has a logarithmic cost. If the interval is logarithmically short, the vanishing forcing cannot create the required positive gain; if it is logarithmically long, the accumulated dissipation forces a loss. This dichotomy replaces any unproved summability assumption on $D_k$.
\end{remark}
\section{The exact fixed-target energy recurrence}\label{sec:fixed}
\subsection{States, coefficients, and the one-step inequality}
Fix once and for all a free parameter $2<\beta<b$. This is possible precisely because $\alpha>3$ and $b=\alpha-1$. For $k\ge1$, define
\begin{equation}\label{eq:coeff-basic}
 n:=k+b,\qquad m:=k+\beta,\qquad
 \ell_k:=\frac{k-1}{\beta},\qquad
 \theta_k:=\frac{k-1}{n}=\frac{\beta\ell_k}{n}.
\end{equation}
The indices $n$ and $m$ always carry the same implicit $k$-dependence; at index $k+1$ they become $n+1$ and $m+1$, respectively. Introduce the states
\begin{equation}\label{eq:states}
 p_k:=x_k-x_{k-1},\qquad \delta_k:=\|p_k\|,\qquad
 u_k:=x_k+\ell_kp_k,\qquad G_k:=\frac{y_k-x_{k+1}}s.
\end{equation}
Thus $x_{k+1}=y_k-sG_k$ and $\theta_kp_k=p_{k+1}+sG_k$. Proximal optimality shows that $G_k$ is the forward--backward mapping associated with $\Phi_k$, whose smooth part $f+\frac{\eps_k}{2}\|\cdot-a\|^2$ has Lipschitz constant
\begin{equation}\label{eq:Lk}
 L_k:=L+\eps_k.
\end{equation}
Choose $k_1$ such that $sL_k<1$ for $k\ge k_1$; this is possible because $sL<1$ and $\eps_k\to0$. The initialization gives $p_1=0$. Since $\theta_1=0$, the increment elimination below is used only for $k\ge2$. All energies at the finitely many earlier indices are finite and play no asymptotic role.

We begin with the forward--backward inequality in a form that retains its exact gradient residual.

\begin{lemma}[One-step inequality]\label{lem:one-step}
For every $v\in\Hil$ and $k\ge1$,
\begin{equation}\label{eq:one-step}
 \Phi_k(x_{k+1})\le\Phi_k(v)+\langle G_k,y_k-v\rangle
 -s\left(1-\frac{sL_k}{2}\right)\|G_k\|^2.
\end{equation}
\end{lemma}
\begin{proof}
Apply the descent lemma at $y_k$ to the smooth part of $\Phi_k$, and apply convexity with the same point to obtain a value at an arbitrary $v$. Since $x_{k+1}-y_k=-sG_k$, this gives
\[
\begin{aligned}
 f(x_{k+1})+\frac{\eps_k}{2}\|x_{k+1}-a\|^2
 &\le f(v)+\frac{\eps_k}{2}\|v-a\|^2\\
 &\quad+\langle\nabla f(y_k)+\eps_k(y_k-a),y_k-v\rangle\\
 &\quad-s\langle\nabla f(y_k)+\eps_k(y_k-a),G_k\rangle
 +\frac{L_ks^2}{2}\|G_k\|^2.
\end{aligned}
\]
Proximal optimality gives
\[
 G_k-\nabla f(y_k)-\eps_k(y_k-a)\in\partial g(x_{k+1}).
\]
Use this subgradient in $g(x_{k+1})\le g(v)+\langle\partial g(x_{k+1}),x_{k+1}-v\rangle$, write $x_{k+1}-v=y_k-v-sG_k$, and add the resulting inequality to the smooth estimate. The gradient and anchor inner products cancel, leaving \eqref{eq:one-step}.
\end{proof}

\subsection{Definition of the fixed-target energy}
For a fixed target $\xi\in\operatorname{dom}g$, introduce the fixed-objective gaps
\begin{equation}\label{eq:fixed-P}
 P_k(\xi):=\Phi_k(x_k)-\Phi_k(\xi),\qquad
 \widehat P_{k+1}(\xi):=\Phi_k(x_{k+1})-\Phi_k(\xi).
\end{equation}
The energy coefficients are
\begin{align}
 A_k&:=\frac{2s}{\beta^2}(m-1)(n-1),\label{eq:A}\\
 \lambda_k&:=-\frac{(b-\beta)(k-1)(k-1+\beta)}{\beta^2n},\label{eq:lambda}\\
 \tau_k&:=\frac{(k-1)(\beta-b)}{\beta n}.\label{eq:tau}
\end{align}
The fixed-target energy and its one-step version, both using the same objective $\Phi_k$ across the recurrence, are
\begin{equation}\label{eq:H-fixed}
 H_k(\xi):=\|u_k-\xi\|^2+A_kP_k(\xi)+\lambda_k\delta_k^2
 -\tau_{k-1}\|x_{k-1}-\xi\|^2
\end{equation}
and
\begin{equation}\label{eq:H-hat}
 \widehat H_{k+1}(\xi):=\|u_{k+1}-\xi\|^2
 +A_{k+1}\widehat P_{k+1}(\xi)+\lambda_{k+1}\delta_{k+1}^2
 -\tau_k\|x_k-\xi\|^2.
\end{equation}
The signs are intentional: $\tau_k<0$ for $k\ge2$, while $\lambda_k$ is negative for all sufficiently large $k$. Consequently, coercivity will require a joint argument rather than a visual inspection of the energy. The auxiliary state $u_k$ replaces the prescribed inertial coefficient by the freely chosen $\ell_k$. The parameter $\beta$ splits the damping budget: $\beta>2$ will produce objective-gap dissipation, while $b-\beta>0$ will produce increment dissipation. The coefficients $\lambda_k$ and $\tau_k$ retain the signed velocity and delayed target-distance terms exactly instead of discarding them prematurely.

Define the dissipation and delayed-remainder coefficients
\begin{align}
 d_k^P&:=\frac{2s}{\beta^2}\left((\beta-2)k+\beta^2-\beta+1-b\right),\label{eq:dP}\\
 d_k^G&:=\frac{s^2mn}{\beta^2}(1-sL_k),\label{eq:dG}\\
 \mu_k&:=\frac{(b-\beta)(k+\beta)(2k+b+1)}{\beta^2(n+1)},\label{eq:mu}\\
 \omega_k&:=\frac{(b+1)(b-\beta)}{\beta n(n-1)}.\label{eq:omega}
\end{align}
All four are positive for sufficiently large $k$; in fact $\mu_k,\omega_k>0$ for $k\ge2$. Their leading orders are
\begin{equation}\label{eq:coefficient-orders}
 A_k\asymp k^2,\quad d_k^P\asymp k,\quad d_k^G\asymp k^2,\quad
 \mu_k\asymp k,\quad \omega_k\asymp k^{-2}.
\end{equation}

\begin{proposition}[Fixed-objective, fixed-target recurrence]\label{prop:fixed-recurrence}
For every $k\ge2$ and every $\xi\in\operatorname{dom}g$, the sequence generated by \eqref{eq:algorithm} satisfies
\begin{equation}\label{eq:fixed-recurrence}
\begin{aligned}
 &\widehat H_{k+1}(\xi)+d_k^PP_k(\xi)
 +d_k^G\|G_k\|^2+\mu_k\delta_{k+1}^2\\
 &\qquad\le H_k(\xi)+\omega_k\|x_{k-1}-\xi\|^2.
\end{aligned}
\end{equation}
\end{proposition}
\begin{proof}
All occurrences of $\xi$ in this proof are fixed. Set $U_k:=u_k-\xi$ and $\kap:=b-\beta>0$.

\medskip
\noindent\emph{Step 1: a weighted objective inequality.}
Apply \eqref{eq:one-step} at $v=x_k$ and $v=\xi$ with weights $1-\beta/n$ and $\beta/n$, respectively. The definitions of $y_k$ and $u_k$ give the exact identity
\begin{equation}\label{eq:weighted-inner}
 \left(1-\frac\beta n\right)(y_k-x_k)+\frac\beta n(y_k-\xi)
 =\frac\beta nU_k.
\end{equation}
Consequently,
\begin{equation}\label{eq:weighted-step}
\begin{aligned}
 \Phi_k(x_{k+1})
 &\le\left(1-\frac\beta n\right)\Phi_k(x_k)+\frac\beta n\Phi_k(\xi)
 +\frac\beta n\langle G_k,U_k\rangle\\
 &\quad-s\left(1-\frac{sL_k}{2}\right)\|G_k\|^2.
\end{aligned}
\end{equation}

\medskip
\noindent\emph{Step 2: eliminate the extrapolated-state inner product.}
Since $1+\ell_{k+1}=m/\beta$ and $p_{k+1}=\theta_kp_k-sG_k$,
\begin{equation}\label{eq:u-recursion-derivation}
\begin{aligned}
 u_{k+1}
 &=x_{k+1}+\ell_{k+1}p_{k+1}\\
 &=x_k+\theta_kp_k-sG_k+\ell_{k+1}(\theta_kp_k-sG_k)\\
 &=x_k+\frac{m\ell_k}{n}p_k-\frac{sm}{\beta}G_k\\
 &=u_k-\eta_kp_k-\frac{sm}{\beta}G_k,
 \qquad \eta_k:=\frac{\kap\ell_k}{n},
\end{aligned}
\end{equation}
where the last equality uses $m\ell_k/n=\ell_k-\eta_k$ because $n-\kap=m$. Thus
\begin{equation}\label{eq:u-recursion}
 u_{k+1}=u_k-\eta_kp_k-\frac{sm}{\beta}G_k.
\end{equation}
Expanding the squared norm of \eqref{eq:u-recursion} gives
\begin{equation}\label{eq:u-norm}
\begin{aligned}
 \|U_{k+1}\|^2={}&\|U_k\|^2+\eta_k^2\delta_k^2
 +\frac{s^2m^2}{\beta^2}\|G_k\|^2-2\eta_k\langle U_k,p_k\rangle\\
 &+\frac{2sm\eta_k}{\beta}\langle p_k,G_k\rangle
 -\frac{2sm}{\beta}\langle U_k,G_k\rangle .
\end{aligned}
\end{equation}
Solving this identity for $\langle G_k,U_k\rangle$ yields
\begin{equation}\label{eq:G-U}
\begin{aligned}
 \frac\beta n\langle G_k,U_k\rangle
 &={}\frac{\beta^2}{2smn}\left(\|U_k\|^2-\|U_{k+1}\|^2\right)
 +\frac{\beta^2\eta_k^2}{2smn}\delta_k^2\\
 &\quad-\frac{\beta^2\eta_k}{smn}\langle U_k,p_k\rangle
 +\frac{\beta\eta_k}{n}\langle p_k,G_k\rangle
 +\frac{sm}{2n}\|G_k\|^2.
\end{aligned}
\end{equation}
The exact increment identity
\begin{equation}\label{eq:increment-identity}
 \delta_{k+1}^2=\theta_k^2\delta_k^2
 -2s\theta_k\langle p_k,G_k\rangle+s^2\|G_k\|^2
\end{equation}
will eliminate the remaining $\langle p_k,G_k\rangle$ term.

Substitute \eqref{eq:G-U} into \eqref{eq:weighted-step}, subtract $\Phi_k(\xi)$, and multiply by $2smn/\beta^2$. Then use \eqref{eq:increment-identity}. More explicitly, the remaining mixed term is replaced by
\begin{equation}\label{eq:mixed-elimination}
 \frac{2sm\eta_k}{\beta}\langle p_k,G_k\rangle
 =\frac{m\kap}{\beta^2}
 \left(\theta_k^2\delta_k^2-\delta_{k+1}^2+s^2\|G_k\|^2\right),
\end{equation}
which follows by multiplying \eqref{eq:increment-identity} by $m\kap/\beta^2$ and using $\theta_k=\beta\ell_k/n$. After this substitution, the coefficient of $\|G_k\|^2$ on the left is
\begin{equation}\label{eq:G-coefficient}
\begin{aligned}
 &\frac{2s^2mn}{\beta^2}\left(1-\frac{sL_k}{2}\right)
 -\frac{s^2m^2}{\beta^2}-\frac{s^2m\kap}{\beta^2}\\
 &\qquad=\frac{s^2mn}{\beta^2}(1-sL_k)=d_k^G,
\end{aligned}
\end{equation}
because $m+\kap=n$. The objective coefficient at $x_k$ becomes $A_k-d_k^P$, because
\begin{equation}\label{eq:A-dP-identity}
 A_k-A_{k+1}\left(1-\frac\beta n\right)=d_k^P,
 \qquad A_{k+1}=\frac{2smn}{\beta^2}.
\end{equation}
The resulting $\|G_k\|^2$ coefficient is exactly $-d_k^G$. After moving the next-iterate terms to the left, one obtains
\begin{equation}\label{eq:pre-delayed}
\begin{aligned}
 &\|U_{k+1}\|^2+A_{k+1}\widehat P_{k+1}(\xi)
 +d_k^PP_k(\xi)+d_k^G\|G_k\|^2
 +\frac{m\kap}{\beta^2}\delta_{k+1}^2\\
 &\le\|U_k\|^2+A_kP_k(\xi)
 +\left(\eta_k^2+\frac{m\kap}{\beta^2}\theta_k^2\right)\delta_k^2
 -2\eta_k\langle U_k,p_k\rangle .
\end{aligned}
\end{equation}

\medskip
\noindent\emph{Step 3: expose the delayed distance.}
Put $w_k:=x_k-\xi$, so that $U_k=w_k+\ell_kp_k$ and $x_{k-1}-\xi=w_k-p_k$. The polarization identity gives
\begin{equation}\label{eq:cross-delayed}
 -2\eta_k\langle U_k,p_k\rangle
 =\eta_k\|w_k-p_k\|^2-\eta_k\|w_k\|^2
 -(\eta_k+2\eta_k\ell_k)\delta_k^2.
\end{equation}
Moreover, $\eta_k=-\tau_k$ and
\begin{equation}\label{eq:lambda-identity}
 \eta_k^2+\frac{m\kap}{\beta^2}\theta_k^2
 -\eta_k-2\eta_k\ell_k
 =-\frac{\kap\ell_k(\ell_k+1)}n=\lambda_k.
\end{equation}
Here the first simplification uses $\kap+m=n$. Substitution of \eqref{eq:cross-delayed}--\eqref{eq:lambda-identity} into \eqref{eq:pre-delayed} gives
\begin{equation}\label{eq:almost-fixed}
\begin{aligned}
 &\widehat H_{k+1}(\xi)+d_k^PP_k(\xi)+d_k^G\|G_k\|^2
 +\left(\frac{m\kap}{\beta^2}-\lambda_{k+1}\right)\delta_{k+1}^2\\
 &\qquad\le H_k(\xi)+(\tau_{k-1}-\tau_k)\|x_{k-1}-\xi\|^2.
\end{aligned}
\end{equation}
Finally, direct substitution of \eqref{eq:lambda} and \eqref{eq:tau} gives
\[
 \frac{m\kap}{\beta^2}-\lambda_{k+1}=\mu_k,
 \qquad \tau_{k-1}-\tau_k=\omega_k.
\]
This proves \eqref{eq:fixed-recurrence}.
\end{proof}

\begin{remark}[The role of the free inertia]
The identity \eqref{eq:A-dP-identity} explains the condition $\beta>2$: the leading term of $d_k^P$ is $2s(\beta-2)k/\beta^2$. The condition $\beta<b$ makes $\kap>0$ and hence $\mu_k>0$. Thus the two sides of $2<\beta<b$ are responsible for objective-gap dissipation and velocity dissipation, respectively.
\end{remark}

\subsection{Joint coercivity and two-index domination}
From now on, $H_k:=H_k(\xi_k)$, with $P_k$ as in \eqref{eq:gap-radius}. This is the actual moving-target energy. The following proposition shows that the signed velocity coefficient causes no asymptotic loss of coercivity.

\begin{proposition}[Joint coercivity]\label{prop:coercivity}
There are $k_0\ge k_1$ and constants $0<c_0\le1\le C_0$ such that, for all $k\ge k_0$,
\begin{equation}\label{eq:coercivity}
\begin{aligned}
 c_0\left(\|x_k-\xi_k\|^2+A_kP_k+\ell_k^2\delta_k^2\right)
 &\le H_k\\
 &\le C_0\left(A_kP_k+\ell_k^2\delta_k^2\right).
\end{aligned}
\end{equation}
In particular, $H_k\ge0$. The threshold $k_0$ and the constants depend only on the fixed parameters $s,L,c,\beta,b$ and the finite initial segment, not on any bound for $(x_k)$.
\end{proposition}
\begin{proof}
Write $z_k:=x_k-\xi_k$ and $\nu_k:=-\tau_{k-1}\ge0$ for large $k$. Splitting $A_kP_k$ into two equal parts and using \eqref{eq:strong-gap} on one half give
\begin{equation}\label{eq:coercive-split}
 H_k\ge\frac12A_kP_k+Q_k(z_k,p_k),
\end{equation}
where, in scalar coefficient notation,
\begin{equation}\label{eq:Q-definition}
\begin{aligned}
 Q_k(z,p):={}&\|z+\ell_kp\|^2+\nu_k\|z-p\|^2+\lambda_k\|p\|^2
 +\frac{A_k\eps_k}{4}\|z\|^2\\
 ={}&a_k\|z\|^2+2c_k\langle z,\ell_kp\rangle+b_k\|\ell_kp\|^2,
\end{aligned}
\end{equation}
with
\[
 a_k:=1+\nu_k+\frac{A_k\eps_k}{4},\qquad
 b_k:=1+\frac{\nu_k+\lambda_k}{\ell_k^2},\qquad
 c_k:=1-\frac{\nu_k}{\ell_k}.
\]
The limits
\begin{equation}\label{eq:matrix-limits}
 a_k\to1+\frac{b-\beta}{\beta}+\frac{sc}{2\beta^2},\qquad b_k\to1,\qquad c_k\to1
\end{equation}
follow from $A_k\sim2sk^2/\beta^2$, $\ell_k\sim k/\beta$, $\lambda_k\sim-(b-\beta)k/\beta^2$, and $\nu_k\to(b-\beta)/\beta$. Hence the symmetric two-dimensional matrix with entries $(a_k,c_k;c_k,b_k)$ is uniformly positive definite for all large $k$: its diagonal entries have positive limits and its determinant tends to $(b-\beta)/\beta+sc/(2\beta^2)>0$. This proves the lower bound in \eqref{eq:coercivity}, including the $\|x_k-\xi_k\|^2$ term.

For the upper bound, discard $\lambda_k\delta_k^2\le0$ and use
\[
 \|z_k+\ell_kp_k\|^2+\nu_k\|z_k-p_k\|^2
 \le(2+2\nu_k)\|z_k\|^2+(2\ell_k^2+2\nu_k)\delta_k^2.
\]
By \eqref{eq:strong-gap},
\[
 \|z_k\|^2\le\frac{2}{\eps_k}P_k
 =\frac{2k^2}{c}P_k\le C A_kP_k
\]
for a parameter-dependent constant $C$, because $A_k$ is of order $k^2$. Since $\nu_k$ is bounded and $\ell_k\to\infty$, the remaining terms are dominated by a constant multiple of $A_kP_k+\ell_k^2\delta_k^2$.
\end{proof}

Define the dissipation block
\begin{equation}\label{eq:D-definition}
 D_k:=A_kP_k+\ell_{k+1}^2\delta_{k+1}^2.
\end{equation}
The upper estimate in \eqref{eq:coercivity}, together with the elementary bounds $A_kP_k\le D_k$ and $\ell_k^2\delta_k^2\le D_{k-1}$, gives the two-index bound required by Lemma~\ref{lem:two-index}:
\begin{equation}\label{eq:H-D-bound}
 H_k\le C_0(D_k+D_{k-1})\qquad(k\ge k_0+1),
\end{equation}
with the same constant $C_0$.
\section{Moving-path recurrence and uniform absorption}\label{sec:moving}
The fixed-target recurrence cannot simply be evaluated at $\xi_k$ on both sides: the target and the regularized objective both change from $k$ to $k+1$. This section computes those two shifts exactly and then absorbs every remainder into the displayed dissipation.

\subsection{The exact moving-path identity}
Recall $h_k=\xi_{k+1}-\xi_k$ and define
\begin{equation}\label{eq:B-definition}
 B_k:=\frac{A_{k+1}}2(\eps_k-\eps_{k+1})
 =\frac{scmn(2k+1)}{\beta^2k^2(k+1)^2}.
\end{equation}
In particular, $0<B_k\le C_q/k$ for a constant $C_q$ depending only on $s,c,\beta,b$.

\begin{proposition}[Exact moving-path estimate]\label{prop:moving-exact}
For every $k\ge2$,
\begin{equation}\label{eq:moving-exact}
\begin{aligned}
 H_{k+1}
 &\le H_k-d_k^PP_k-d_k^G\|G_k\|^2-\mu_k\delta_{k+1}^2
 +B_kq_{k+1,+}\\
 &\quad+\omega_k\|x_{k-1}-\xi_k-h_k\|^2
 -2\langle u_k-\xi_k-\tau_{k-1}(x_{k-1}-\xi_k),h_k\rangle\\
 &\quad+(1-\tau_{k-1})\|h_k\|^2
 -(A_k-d_k^P)\bigl(\Phi_k(\xi_{k+1})-\Phi_k(\xi_k)\bigr).
\end{aligned}
\end{equation}
The last term on the right is nonpositive.
\end{proposition}
\begin{proof}
Write $\widetilde P_k:=P_k(\xi_{k+1})$ and apply \eqref{eq:fixed-recurrence} at $\xi=\xi_{k+1}$.

\medskip
\noindent\emph{Step 1: shift the next objective.}
The fixed next energy differs from $H_{k+1}$ only in the objective term. Since
\[
 \widehat P_{k+1}(\xi_{k+1})-P_{k+1}
 =\frac{\eps_k-\eps_{k+1}}2
 \left(\|x_{k+1}-a\|^2-R_{k+1}^2\right),
\]
the definition $q_{k+1}=R_{k+1}^2-\|x_{k+1}-a\|^2$ gives
\begin{equation}\label{eq:next-objective-shift}
 H_{k+1}=\widehat H_{k+1}(\xi_{k+1})+B_kq_{k+1}
 \le\widehat H_{k+1}(\xi_{k+1})+B_kq_{k+1,+}.
\end{equation}

\medskip
\noindent\emph{Step 2: shift the current target.}
Put $U_k:=u_k-\xi_k$, $v_k:=x_{k-1}-\xi_k$, and
$\Delta\Phi_k:=\Phi_k(\xi_{k+1})-\Phi_k(\xi_k)\ge0$. Then
$\widetilde P_k=P_k-\Delta\Phi_k$, and direct expansion gives
\begin{equation}\label{eq:H-shift}
\begin{aligned}
 H_k(\xi_{k+1})
 ={}&H_k-2\langle U_k-\tau_{k-1}v_k,h_k\rangle
 +(1-\tau_{k-1})\|h_k\|^2-A_k\Delta\Phi_k.
\end{aligned}
\end{equation}
Combining \eqref{eq:fixed-recurrence}, \eqref{eq:next-objective-shift}, and \eqref{eq:H-shift} yields \eqref{eq:moving-exact}, because
\begin{equation}\label{eq:objective-cancellation}
 -A_k\Delta\Phi_k-d_k^P\widetilde P_k
 =-d_k^PP_k-(A_k-d_k^P)\Delta\Phi_k.
\end{equation}
Finally, \eqref{eq:A-dP-identity} shows that
$A_k-d_k^P=A_{k+1}(1-\beta/n)>0$ for $k\ge2$, so the last term is nonpositive.
\end{proof}

The cancellation \eqref{eq:objective-cancellation} is essential. It removes the moving-target gap $\widetilde P_k$ completely and leaves a nonpositive multiple of the nonnegative path-objective increment $\Delta\Phi_k$. No growth estimate for $F$ along the path is required.

\subsection{Absorption into the dissipation block}
\begin{proposition}[Uniform absorption]\label{prop:absorption}
There are $k_*\ge k_0+1$, $\rho>0$, $C_q>0$, and a nonnegative sequence $(\zeta_k)$ such that $\sum_k\zeta_k<\infty$ and, for all $k\ge k_*$,
\begin{equation}\label{eq:final-recurrence}
 H_{k+1}\le H_k-\frac\rho kD_k+\frac{C_q}{k}q_{k+1,+}+\zeta_k.
\end{equation}
One may take $\zeta_k=C_\gamma k\|h_k\|^2$ after fixing sufficiently small absorption parameters $\gamma_1,\gamma_2>0$. No multiplicative $(1+\sigma_k)H_k$ term is needed; equivalently, $\sigma_k=0$.
\end{proposition}
\begin{proof}
Discard the nonpositive last term of \eqref{eq:moving-exact}. By \eqref{eq:coefficient-orders}, there are constants $d_0,g_0,m_0>0$ such that, for all large $k$,
\begin{equation}\label{eq:coefficient-lower}
 d_k^P\ge d_0k,\qquad d_k^G\ge g_0k^2,\qquad \mu_k\ge m_0k.
\end{equation}
Here $\beta>2$, $sL_k<1$ eventually, and $\beta<b$ are the three relevant parameter conditions.

Set $z_k:=x_k-\xi_k$, $U_k:=z_k+\ell_kp_k$, and $v_k:=z_k-p_k$. The path inner product in \eqref{eq:moving-exact} expands as
\[
 -2(1-\tau_{k-1})\langle z_k,h_k\rangle
 -2(\ell_k+\tau_{k-1})\langle p_k,h_k\rangle .
\]
Since $\tau_{k-1}$ is bounded, Young's inequality with an arbitrary $\gamma_1>0$ and the strong-gap estimate \eqref{eq:strong-gap} give
\begin{equation}\label{eq:e-absorption}
\begin{aligned}
 2|1-\tau_{k-1}|\,|\langle z_k,h_k\rangle|
 &\le C\left(\frac{\gamma_1}{k}\|z_k\|^2+\frac{k}{\gamma_1}\|h_k\|^2\right)\\
 &\le \frac{2C\gamma_1}{c}kP_k+\frac{Ck}{\gamma_1}\|h_k\|^2.
\end{aligned}
\end{equation}
The constant $C$ is independent of $\gamma_1$ and $k$. Since $\ell_k+\tau_{k-1}=O(k)$, similarly, for every $\gamma_2>0$,
\begin{equation}\label{eq:p-absorption}
 2|\ell_k+\tau_{k-1}|\,|\langle p_k,h_k\rangle|
 \le C\gamma_2 k\delta_k^2+\frac{Ck}{\gamma_2}\|h_k\|^2.
\end{equation}
The increment identity $\theta_kp_k=p_{k+1}+sG_k$ and $\theta_k\ge1/2$ for large $k$ imply
\begin{equation}\label{eq:p-control}
 \delta_k^2\le C(\delta_{k+1}^2+\|G_k\|^2).
\end{equation}
Thus the first term on the right of \eqref{eq:p-absorption} is bounded by
$C\gamma_2k\delta_{k+1}^2+C\gamma_2k\|G_k\|^2$.

The delayed remainder satisfies
\begin{equation}\label{eq:delayed-absorption}
\begin{aligned}
 \omega_k\|v_k-h_k\|^2
 &\le2\omega_k\|v_k\|^2+2\omega_k\|h_k\|^2\\
 &\le4\omega_k\|z_k\|^2+4\omega_k\delta_k^2+2\omega_k\|h_k\|^2\\
 &\le C(P_k+k^{-2}\delta_{k+1}^2+k^{-2}\|G_k\|^2)+Ck\|h_k\|^2,
\end{aligned}
\end{equation}
where $\omega_k=O(k^{-2})$, $\|z_k\|^2\le2k^2P_k/c$, and \eqref{eq:p-control} were used. Finally,
\[
 (1-\tau_{k-1})\|h_k\|^2\le Ck\|h_k\|^2.
\]

Now choose $\gamma_1$ small enough that the $kP_k$ term in \eqref{eq:e-absorption} is at most $\frac14d_0kP_k$, and choose $\gamma_2$ small enough that the velocity and gradient remainders proportional to $\gamma_2$ are at most one quarter of their respective dissipation terms. After that choice, enlarge $k$ so that the $\gamma$-independent remainders $CP_k$, $Ck^{-2}\delta_{k+1}^2$, and $Ck^{-2}\|G_k\|^2$ are likewise dominated. In particular, the gradient remainder of order $k\|G_k\|^2$ is dominated by $d_k^G\asymp k^2$. This yields
\[
 H_{k+1}\le H_k-\frac{d_0}{2}kP_k-\frac{m_0}{2}k\delta_{k+1}^2
 +\frac{C_q}{k}q_{k+1,+}+\zeta_k,
\]
with $\zeta_k=C_\gamma k\|h_k\|^2$. Lemma~\ref{lem:path-summability} proves $\sum_k\zeta_k<\infty$. Since $A_k\asymp k^2$ and $\ell_{k+1}\asymp k$, there is $\rho>0$ such that
\[
 \frac{d_0}{2}kP_k+\frac{m_0}{2}k\delta_{k+1}^2
 \ge\frac\rho k(A_kP_k+\ell_{k+1}^2\delta_{k+1}^2)=\frac\rho kD_k.
\]
This proves \eqref{eq:final-recurrence}.
\end{proof}

\begin{remark}[Uniformity of the constants]
The absorption thresholds depend on $s,c,\beta,b$, the Lipschitz constant $L$, and fixed coefficient bounds. They do not depend on $\sup_k\|x_k\|$, on $\sup_kP_k$, or on a bound for the path increments. All terms involving $h_k$ are collected in the single summable sequence $\zeta_k$; none is hidden in a multiplicative energy factor.
\end{remark}
\section{Proof of the main theorem}\label{sec:proof}
We prove Theorem~\ref{thm:main} in three stages. The first uses only the uniform bound $q_{k,+}\le R_*^2$ and therefore does not assume that $(x_k)$ is bounded. The second derives objective convergence and upgrades the radial forcing. The third applies the two-index logarithmic lemma and extracts all stated conclusions.

\subsection{Stage 1: the scaled energy tends to zero}
For $k\ge k_*$, drop the nonnegative dissipation in \eqref{eq:final-recurrence} and use \eqref{eq:q-positive-bound} to obtain
\[
 H_{k+1}\le H_k+\frac{C_qR_*^2}{k}+\zeta_k.
\]
Set $E_k:=H_k/k^2$. Proposition~\ref{prop:coercivity} gives $E_k\ge0$. Division by $(k+1)^2$ yields
\begin{equation}\label{eq:scaled-bootstrap}
 E_{k+1}\le\frac{k^2}{(k+1)^2}E_k
 +\frac{C_qR_*^2}{k(k+1)^2}+\frac{\zeta_k}{(k+1)^2}.
\end{equation}
The damping factor satisfies
\[
 \frac{k^2}{(k+1)^2}
 =1-\frac2k+\frac{3k+2}{k(k+1)^2},
\]
so its correction to $1-2/k$ is nonnegative and summable. Both forcing sequences on the right of \eqref{eq:scaled-bootstrap} are summable: the first is $O(k^{-3})$, and the second is dominated termwise by $\zeta_k$. Lemma~\ref{lem:damped} with $\alpha=2$ therefore gives
\begin{equation}\label{eq:E-to-zero}
 \frac{H_k}{k^2}=E_k\longrightarrow0.
\end{equation}
Because $A_k\asymp k^2$, the lower bound in Proposition~\ref{prop:coercivity} now implies
\begin{equation}\label{eq:P-to-zero}
 P_k\longrightarrow0.
\end{equation}
At this point we have regularized-gap convergence, but not yet strong convergence or even boundedness of $(x_k)$.

\subsection{Stage 2: objective convergence and vanishing radial forcing}
Comparison of $\Phi_k(\xi_k)$ with $\Phi_k(z)$ gives
\begin{equation}\label{eq:objective-P-lower}
\begin{aligned}
 P_k
 &\ge\Phi_k(x_k)-\Phi_k(z)\\
 &=F(x_k)-F^*+\frac{\eps_k}{2}
 \left(\|x_k-a\|^2-R_*^2\right).
\end{aligned}
\end{equation}
Consequently,
\begin{equation}\label{eq:objective-basic-bound}
 0\le F(x_k)-F^*
 \le P_k+\frac{\eps_k}{2}\bigl(R_*^2-\|x_k-a\|^2\bigr)
 \le P_k+\frac{\eps_k}{2}R_*^2.
\end{equation}
By \eqref{eq:P-to-zero} and $\eps_k\to0$, $F(x_k)\to F^*$. Lemma~\ref{lem:negative-radial} and \eqref{eq:q-positive-bound} therefore give
\begin{equation}\label{eq:q-to-zero}
 q_{k,+}\longrightarrow0,
 \qquad \frac{C_q}{k}q_{k+1,+}=o(k^{-1}).
\end{equation}
No boundedness or weak-convergence property of the full sequence $(x_k)$ was used. The negative-radial argument extracts a bounded subsequence only after objective convergence has been established.

\subsection{Stage 3: the two-index limit and the final rates}
The recurrence \eqref{eq:final-recurrence}, the domination \eqref{eq:H-D-bound}, and \eqref{eq:q-to-zero} satisfy all hypotheses of Lemma~\ref{lem:two-index}, with
\[
 r_k=\frac{C_q}{k}q_{k+1,+}=o(k^{-1}).
\]
Hence
\begin{equation}\label{eq:H-to-zero}
 H_k\longrightarrow0.
\end{equation}
The lower bound in Proposition~\ref{prop:coercivity} then gives the three simultaneous limits
\begin{equation}\label{eq:three-limits}
 \|x_k-\xi_k\|\to0,\qquad A_kP_k\to0,\qquad \ell_k^2\delta_k^2\to0.
\end{equation}
Lemma~\ref{lem:path-convergence} and the first limit prove
\[
 x_k\longrightarrow z=P_Sa
\]
in the strong topology. Since $A_k/k^2\to2s/\beta^2$ and $\ell_k/k\to1/\beta$, the other two limits in \eqref{eq:three-limits} give
\begin{equation}\label{eq:P-delta-little-o}
 k^2P_k\to0,
 \qquad k\delta_k\to0.
\end{equation}
The second assertion is exactly $\|x_{k+1}-x_k\|=o(k^{-1})$ after an index shift.

Finally, strong convergence $x_k\to z$ implies $\|x_k-a\|^2-R_*^2\to0$. Multiplying \eqref{eq:objective-P-lower} by $k^2$ yields
\[
 0\le k^2(F(x_k)-F^*)
 \le k^2P_k-\frac c2\left(\|x_k-a\|^2-R_*^2\right)\longrightarrow0.
\]
Thus $F(x_k)-F^*=o(k^{-2})$, completing the proof of Theorem~\ref{thm:main}. \qed

\begin{corollary}[Convergence of the extrapolated states]\label{cor:states}
Under the hypotheses of Theorem~\ref{thm:main},
\[
 y_k\longrightarrow P_Sa,
 \qquad u_k\longrightarrow P_Sa,
 \qquad \|y_k-x_k\|=o(k^{-1}).
\]
\end{corollary}
\begin{proof}
The identity $y_k-x_k=\theta_kp_k$ and $0<\theta_k<1$ give
\[
 \|y_k-x_k\|\le\delta_k=o(k^{-1}).
\]
Also $u_k-x_k=\ell_kp_k$ and $\ell_k^2\delta_k^2\to0$ by \eqref{eq:three-limits}. Both convergence statements follow from $x_k\to P_Sa$.
\end{proof}

\begin{remark}[Where the critical exponent is used]
Weighted path summability is not, by itself, specific to the quadratic Tikhonov exponent. Indeed, if $\eps_k=ck^{-p}$ with $p>0$, then
\[
 \frac{\eps_k-\eps_{k+1}}{\eps_k+\eps_{k+1}}\sim\frac{p}{2k},
\]
so \eqref{eq:path-one-step} gives $k\|h_k\|^2\le C_p(R_{k+1}^2-R_k^2)$ for all sufficiently large $k$, and hence $\sum_k k\|h_k\|^2<\infty$. The quadratic exponent is used in the compatibility of this path estimate with the remaining energy estimates. At $p=2$, $A_k\asymp k^2$, $A_k\eps_k\asymp1$, and $B_k=O(k^{-1})$; these relations yield the uniform upper coercivity estimate and the objective-parameter shift $O(k^{-1})q_{k+1,+}$. The initial bounded forcing $q_{k,+}\le R_*^2$ is sufficient for the $k^{-2}$-scaled bootstrap; the later improvement $q_{k,+}\to0$ is exactly what raises that forcing to the logarithmically effective $o(k^{-1})$ scale. Thus the argument neither treats the anchor as an ordinary summable perturbation on the accelerated energy scale nor assumes a one-sided geometric configuration.
\end{remark}

\begin{remark}[Parameter range and anchor strength]
The only role of $\alpha>3$ is to permit $2<\beta<b=\alpha-1$. This gives $d_k^P\asymp(\beta-2)k$, while $\beta<b$ gives $\mu_k\asymp2(b-\beta)k/\beta^2$. The anchor coefficient $c>0$ enters the coercivity determinant only through the nonnegative contribution $sc/(2\beta^2)$ to its limit $(b-\beta)/\beta+sc/(2\beta^2)$. No lower bound on $c$ or upper bound on $sc$ is used; the role of the condition $c>0$ is to define the Tikhonov path on which the selection mechanism rests. The stepsize restriction enters through $sL_k<1$ eventually, which follows from $sL<1$ and $\eps_k\to0$.
\end{remark}

\section{Stability under additive errors}\label{sec:errors}
Let $(e_k)_{k\ge1}\subset\Hil$ be an arbitrary sequence of computational errors. We consider the inexact iteration
\begin{equation}\label{eq:inexact-algorithm}
\boxed{
\begin{aligned}
 y_k&=x_k+\frac{k-1}{k+b}(x_k-x_{k-1}),\qquad b=\alpha-1,\\
 x_{k+1}&=\prox_{sg}\!\left(y_k-s\nabla f(y_k)-\frac{sc}{k^2}(y_k-a)+e_k\right),
\end{aligned}}
\end{equation}
with $x_0=x_1\in\operatorname{dom}g$. Thus $e_k$ is the total additive error in the argument of the proximal operator. If $\widetilde{\nabla f}_k$ denotes the approximation to $\nabla f(y_k)$ used in the forward step, then the corresponding total additive error is
\[
 e_k=s\bigl(\nabla f(y_k)-\widetilde{\nabla f}_k\bigr).
\]
Writing $\chi_k:=\nabla f(y_k)-\widetilde{\nabla f}_k$ for the gradient error, this reads $e_k=s\chi_k$.

\begin{theorem}[Strong selection and rates under additive errors]\label{thm:inexact}
Let $\alpha>3$, $c>0$, $a\in\Hil$, $0<s<1/L$, and suppose that
\begin{equation}\label{eq:error-tolerance}
 \sum_{k=1}^{\infty}k\|e_k\|<\infty.
\end{equation}
For any $2<\beta<b=\alpha-1$, define the auxiliary state
\[
 u_k:=x_k+\frac{k-1}{\beta}(x_k-x_{k-1}).
\]
Then the sequence generated by \eqref{eq:inexact-algorithm} satisfies
\[
 x_k\longrightarrow P_Sa\quad\text{strongly},\qquad
 F(x_k)-F^*=o(k^{-2}),\qquad
 \|x_{k+1}-x_k\|=o(k^{-1}).
\]
Moreover, $y_k\to P_Sa$, $u_k\to P_Sa$, and $\|y_k-x_k\|=o(k^{-1})$.
\end{theorem}

\begin{proof}
We retain all the notation introduced in Sections~\ref{sec:path}--\ref{sec:proof}. The Tikhonov path and Lemmas~\ref{lem:path-convergence}--\ref{lem:negative-radial} are independent of the iterates. Proposition~\ref{prop:coercivity} is a pointwise statement about the energy evaluated at the current states, so it also remains valid for the states produced by \eqref{eq:inexact-algorithm}. Constants $C$ below have the same parameter dependence as before and never depend on a trajectory bound.

\medskip
\noindent\emph{Step 1: the perturbed one-step and fixed-target inequalities.}
Define $G_k:=(y_k-x_{k+1})/s$, so that $x_{k+1}=y_k-sG_k$ and $\theta_kp_k=p_{k+1}+sG_k$ still hold. Proximal optimality in \eqref{eq:inexact-algorithm} gives
\begin{equation}\label{eq:inexact-prox-opt}
 G_k-\nabla f(y_k)-\eps_k(y_k-a)+\frac{e_k}{s}\in\partial g(x_{k+1}).
\end{equation}
Using \eqref{eq:inexact-prox-opt} in place of the exact proximal optimality relation in the proof of Lemma~\ref{lem:one-step}, and writing $x_{k+1}-v=y_k-v-sG_k$, we obtain, for every $v\in\Hil$ and $k\ge1$,
\begin{equation}\label{eq:inexact-one-step}
\begin{aligned}
 \Phi_k(x_{k+1})\le{}&\Phi_k(v)+\langle G_k,y_k-v\rangle
 -s\left(1-\frac{sL_k}{2}\right)\|G_k\|^2\\
 &+\left\langle\frac{e_k}{s},x_{k+1}-v\right\rangle .
\end{aligned}
\end{equation}
Indeed, the smooth descent and convexity estimates are unchanged, while the additional subgradient term is exactly $\langle e_k/s,x_{k+1}-v\rangle$.

Fix $\xi\in\operatorname{dom}g$ and put $U_k:=u_k-\xi$. Apply \eqref{eq:inexact-one-step} at $v=x_k$ and $v=\xi$ with weights $1-\beta/n$ and $\beta/n$, respectively. The weighted error term is
\[
\begin{aligned}
 &\left\langle\frac{e_k}{s},x_{k+1}-
 \left(\left(1-\frac\beta n\right)x_k+\frac\beta n\xi\right)\right\rangle\\
 &\qquad=\left\langle\frac{e_k}{s},\frac\beta nU_k-sG_k\right\rangle
 =\frac{\beta}{sn}\langle e_k,U_k\rangle-\langle e_k,G_k\rangle ,
\end{aligned}
\]
where \eqref{eq:weighted-inner} was used. Since $A_{k+1}\beta/(sn)=2m/\beta$, repeating Steps 1--3 of the proof of Proposition~\ref{prop:fixed-recurrence} gives, for every $k\ge2$ and every fixed $\xi\in\operatorname{dom}g$,
\begin{equation}\label{eq:inexact-fixed-recurrence}
\begin{aligned}
 &\widehat H_{k+1}(\xi)+d_k^PP_k(\xi)+d_k^G\|G_k\|^2+\mu_k\delta_{k+1}^2\\
 &\qquad\le H_k(\xi)+\omega_k\|x_{k-1}-\xi\|^2
 +\frac{2m}{\beta}\langle e_k,U_k\rangle-A_{k+1}\langle e_k,G_k\rangle .
\end{aligned}
\end{equation}
All the exact algebraic identities are unchanged because they depend only on $x_{k+1}=y_k-sG_k$.

\medskip
\noindent\emph{Step 2: the moving-path recurrence with errors.}
Apply \eqref{eq:inexact-fixed-recurrence} at $\xi=\xi_{k+1}$. Let $\mathcal R_k$ denote the entire right-hand side of \eqref{eq:moving-exact}, and now set $U_k:=u_k-\xi_k$. The objective and target shifts in the proof of Proposition~\ref{prop:moving-exact} are unaffected by $e_k$. Since $u_k-\xi_{k+1}=U_k-h_k$, we obtain
\begin{equation}\label{eq:inexact-moving-recurrence}
 H_{k+1}\le\mathcal R_k
 +\frac{2m}{\beta}\langle e_k,U_k-h_k\rangle
 -A_{k+1}\langle e_k,G_k\rangle .
\end{equation}
Thus the errors add only the three displayed inner products to the exact moving-path estimate.

\medskip
\noindent\emph{Step 3: absorption and summability of the error remainders.}
For $k\ge k_0$, the lower bound in Proposition~\ref{prop:coercivity} gives
\[
 \|U_k\|^2
 \le2\|x_k-\xi_k\|^2+2\ell_k^2\delta_k^2
 \le CH_k.
\]
Consequently, after enlarging $C$ and the starting index,
\begin{equation}\label{eq:error-U-bound}
 \frac{2m}{\beta}|\langle e_k,U_k\rangle|
 \le Ck\|e_k\|\sqrt{H_k}
 \le \sigma_k(1+H_k),
 \qquad \sigma_k:=Ck\|e_k\|.
\end{equation}
Also,
\begin{equation}\label{eq:error-h-bound}
 \frac{2m}{\beta}|\langle e_k,h_k\rangle|
 \le Ck\|e_k\|\|h_k\|
 \le C\left(k^2\|e_k\|^2+\|h_k\|^2\right).
\end{equation}
For the gradient term, Young's inequality and $A_{k+1}^2/d_k^G\le Ck^2$ for large $k$ give
\begin{equation}\label{eq:error-G-bound}
 A_{k+1}|\langle e_k,G_k\rangle|
 \le\frac14d_k^G\|G_k\|^2+\frac{A_{k+1}^2}{d_k^G}\|e_k\|^2
 \le\frac14d_k^G\|G_k\|^2+Ck^2\|e_k\|^2.
\end{equation}
The absorption in Proposition~\ref{prop:absorption} leaves at least one half of $d_k^G\|G_k\|^2$ unused. Therefore \eqref{eq:error-G-bound} can be absorbed without changing the choices of $\gamma_1,\gamma_2$ or the remaining objective and velocity dissipation. Substituting \eqref{eq:error-U-bound}--\eqref{eq:error-G-bound} into \eqref{eq:inexact-moving-recurrence} and repeating the proof of Proposition~\ref{prop:absorption} yield $k_e\ge k_*$ such that, for all $k\ge k_e$,
\begin{equation}\label{eq:inexact-final-recurrence}
 H_{k+1}\le(1+\sigma_k)H_k-\frac\rho kD_k
 +\frac{C_q}{k}q_{k+1,+}+\widehat\zeta_k,
\end{equation}
where
\[
 \widehat\zeta_k:=\zeta_k+
 C\left(k\|e_k\|+k^2\|e_k\|^2+\|h_k\|^2\right).
\]
By \eqref{eq:error-tolerance}, $\sum_k\sigma_k<\infty$. Moreover, $M:=\sup_k k\|e_k\|<\infty$, and hence
\[
 \sum_k k^2\|e_k\|^2
 \le M\sum_k k\|e_k\|<\infty.
\]
Since also $\sum_k\|h_k\|^2\le\sum_k k\|h_k\|^2<\infty$, Lemma~\ref{lem:path-summability} and the summability of $(\zeta_k)$ imply $\sum_k\widehat\zeta_k<\infty$. The domination \eqref{eq:H-D-bound} is unchanged.

\medskip
\noindent\emph{Step 4: the scaled-energy bootstrap with errors.}
Dropping the dissipation in \eqref{eq:inexact-final-recurrence} and using $q_{k+1,+}\le R_*^2$ give
\[
 H_{k+1}\le(1+\sigma_k)H_k+\frac{C_qR_*^2}{k}+\widehat\zeta_k.
\]
For $E_k:=H_k/k^2$, division by $(k+1)^2$ yields
\begin{equation}\label{eq:inexact-scaled-bootstrap}
 E_{k+1}\le\left(1-\frac2k+\widetilde\sigma_k\right)E_k
 +\frac{C_qR_*^2}{k(k+1)^2}+\frac{\widehat\zeta_k}{(k+1)^2},
\end{equation}
where
\[
 \widetilde\sigma_k:=\frac{3k+2}{k(k+1)^2}
 +\frac{k^2}{(k+1)^2}\sigma_k
\]
is nonnegative and summable. Both additive forcings in \eqref{eq:inexact-scaled-bootstrap} are summable. Lemma~\ref{lem:damped} therefore gives $E_k\to0$, and the lower coercivity bound gives $P_k\to0$. No boundedness assumption on $(x_k)$ was used.

\medskip
\noindent\emph{Step 5: objective convergence and the radial upgrade.}
The comparison \eqref{eq:objective-P-lower} and the estimate \eqref{eq:objective-basic-bound} depend only on $P_k$ and $\eps_k$. Hence $F(x_k)\to F^*$. Lemma~\ref{lem:negative-radial} and \eqref{eq:q-positive-bound} then give $q_{k,+}\to0$, so
\[
 r_k:=\frac{C_q}{k}q_{k+1,+}=o(k^{-1}).
\]

\medskip
\noindent\emph{Step 6: removal of the multiplicative factor and conclusion.}
Set
\[
 \Lambda_k:=\prod_{j=k}^{\infty}(1+\sigma_j),\qquad k\ge k_e.
\]
Since $\sum_k\sigma_k<\infty$, the product converges and there is $\Lambda<\infty$ such that $1\le\Lambda_k\le\Lambda$ and $\Lambda_k\to1$. Moreover, $\Lambda_k=(1+\sigma_k)\Lambda_{k+1}$. Put $\bar H_k:=\Lambda_kH_k$. Multiplying \eqref{eq:inexact-final-recurrence} by $\Lambda_{k+1}$ and using $1\le\Lambda_{k+1}\le\Lambda$ give
\begin{equation}\label{eq:inexact-rescaled-recurrence}
 \bar H_{k+1}\le\bar H_k-\frac\rho kD_k
 +\Lambda r_k+\Lambda\widehat\zeta_k .
\end{equation}
Here $\Lambda r_k=o(k^{-1})$, $\sum_k\Lambda\widehat\zeta_k<\infty$, and \eqref{eq:H-D-bound} implies
\[
 \bar H_k\le\Lambda C_0(D_k+D_{k-1}).
\]
Lemma~\ref{lem:two-index} therefore yields $\bar H_k\to0$. Since $0\le H_k\le\bar H_k$, also $H_k\to0$. The remainder of the proof of Theorem~\ref{thm:main}, starting with \eqref{eq:three-limits}, now applies verbatim and proves all three conclusions. The additional assertions for $y_k$ and $u_k$ follow exactly as in Corollary~\ref{cor:states}.
\end{proof}

\begin{remark}[The error tolerance]
Writing $e_k=s\chi_k$ for a gradient error $\chi_k$ turns \eqref{eq:error-tolerance} into $\sum_k k\|\chi_k\|<\infty$, the tolerance used for inexact Nesterov-type forward--backward methods in \cite{AttouchPeypouquet2016}. The argument above preserves not only the two little-$o$ rates but also strong hierarchical selection to $P_Sa$. At no point does it require boundedness of the inexact trajectory.
\end{remark}

\section{Numerical illustrations}\label{sec:numerics}
We illustrate the selection statement on two deterministic examples with nonunique minimizers. Both use $\alpha=4$, $s=0.95$, and $100{,}000$ iterations. The experiments display the anchor-dependent limits and the decay of the scaled objective and velocity residuals. All values and figures were generated directly from \eqref{eq:algorithm}, without random sampling or external data.

\subsection{Slow finite-horizon selection for a small anchor coefficient}
\begin{figure}[htbp]
\centering
\includegraphics[width=0.96\linewidth]{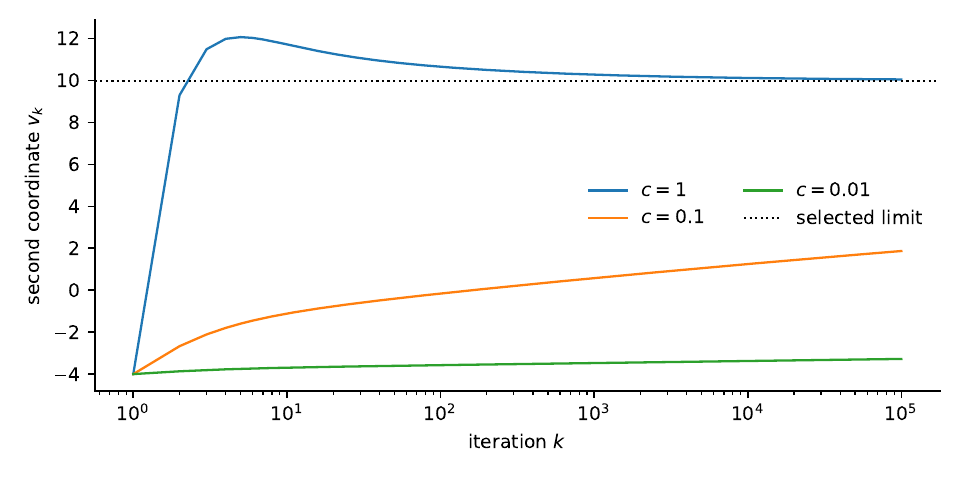}
\caption{Example with $S=\{0\}\times\R$: second coordinate $v_k$ for three anchor coefficients. The dotted line is the second coordinate of $P_Sa=(0,10)$. Small $c$ can produce a long selection transient even though the asymptotic theorem applies for every $c>0$.}
\label{fig:slow-selection}
\end{figure}

Let $\Hil=\R^2$, $f=0$, and $g(u,v)=\iota_{\{u=0\}}$, the indicator of the line $u=0$. Then $S=\{0\}\times\R$, and for the anchor $a=(0,10)$ the hierarchical limit is $P_Sa=(0,10)$. We take $x_0=x_1=(0,-4)$, which belongs to $\operatorname{dom}g$, and $c\in\{1,0.1,0.01\}$. Since $\nabla f=0$ and $\prox_{sg}(u,v)=(0,v)$, the second coordinate evolves by
\[
 v_{k+1}=y_{k,2}-\frac{sc}{k^2}(y_{k,2}-10),
\]
while the first coordinate is reset to zero. Figure~\ref{fig:slow-selection} shows that $c=1$ is close to the selected limit at the terminal iteration, whereas the two smaller coefficients still exhibit substantial finite-horizon transients. Table~\ref{tab:terminal-results} reports the terminal second coordinates and scaled velocities for all three coefficients.

\subsection{Opposite endpoint selection and scaled residuals}
\begin{figure}[htbp]
\centering
\includegraphics[width=\linewidth]{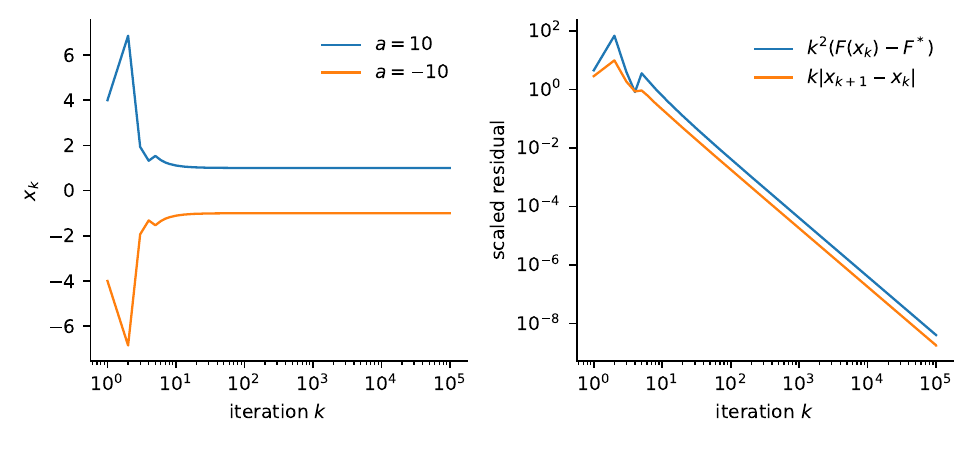}
\caption{Example with $S=[-1,1]$: anchor-dependent endpoint selection (left), and scaled objective and velocity residuals for $a=10$ (right).}
\label{fig:selection-rates}
\end{figure}

Let $\Hil=\R$, $g=0$, and
\[
 f(x)=\frac12\operatorname{dist}(x,[-1,1])^2,
\]
so that $L=1$ and $S=[-1,1]$. With $c=1$ and $x_0=x_1=4\operatorname{sign}(a)$, the anchors $a=10$ and $a=-10$ select the opposite endpoints $1$ and $-1$, respectively. Figure~\ref{fig:selection-rates} displays both trajectories and, for $a=10$, the scaled objective gap $k^2(F(x_k)-F^*)$ and scaled velocity $k|x_{k+1}-x_k|$. Table~\ref{tab:terminal-results} gives the terminal states and both scaled residuals for the two anchors. The plotted decay illustrates the two little-$o$ conclusions of Theorem~\ref{thm:main} over the computed horizon.

\begin{table}[htbp]
\centering
\caption{Terminal results at $k=100{,}000$. The line example reports the second coordinate $v_k$ and has zero objective gap at every iterate. The endpoint example uses $c=1$; its two anchors produce identical scaled residuals by symmetry.}
\label{tab:terminal-results}
\small
\setlength{\tabcolsep}{3pt}
\begin{tabularx}{\linewidth}{@{}l c c c >{\centering\arraybackslash}X >{\centering\arraybackslash}X >{\centering\arraybackslash}X@{}}
\toprule
Example & $c$ & $a$ & $P_Sa$ & Terminal state & $k^2(F(x_k)-F^*)$ & $k|x_{k+1}-x_k|$ \\
\midrule
Line & $1$ & $(0,10)$ & $(0,10)$ & $10.053625676$ & $0$ & $1.929696083\times10^{-2}$ \\
Line & $0.1$ & $(0,10)$ & $(0,10)$ & $1.870051745$ & $0$ & $2.602366876\times10^{-1}$ \\
Line & $0.01$ & $(0,10)$ & $(0,10)$ & $-3.277367750$ & $0$ & $4.209137114\times10^{-2}$ \\
\addlinespace
Endpoint & $1$ & $10$ & $1$ & $1.000000000900$ & $4.050162542\times10^{-9}$ & $1.798561300\times10^{-9}$ \\
Endpoint & $1$ & $-10$ & $-1$ & $-1.000000000900$ & $4.050162542\times10^{-9}$ & $1.798561300\times10^{-9}$ \\
\bottomrule
\end{tabularx}
\end{table}

\section{Conclusion}\label{sec:conclusion}
For the standard single-anchor inertial forward--backward iteration at the critical Tikhonov exponent, we have proved strong last-iterate selection with objective rate $o(k^{-2})$ and discrete-velocity rate $o(k^{-1})$. The proof isolates the moving-path and radial-forcing obstructions at this exponent and resolves them through weighted path summability, a noncircular radial bootstrap, joint energy coercivity, and a two-index logarithmic limit principle. The same mechanism preserves selection and both rates under additive errors with $\sum_k k\|e_k\|<\infty$.

\section*{Statements and declarations}
\noindent\textbf{Funding.} No funding was received to assist with the preparation of this manuscript.

\smallskip
\noindent\textbf{Competing interests.} The author has no relevant financial or non-financial interests to disclose.

\smallskip
\noindent\textbf{Use of generative artificial intelligence.} During the revision of this manuscript, the author used AI tool for language editing, restructuring of the introduction and conclusion, LaTeX format conversion, proof-checking diagnostics, and assistance with reconstructing and checking the numerical reproduction code. All AI-assisted output is subject to the author's review, and the author takes full responsibility for the final content of the manuscript.

\smallskip
\noindent\textbf{Data and reproducibility.} No external datasets were used. The deterministic source code accompanying this manuscript regenerates every reported numerical value and both figures from the stated parameters.

\end{document}